\documentclass[10pt]{amsart}
\usepackage[pagebackref]{hyperref}

\renewcommand{\d}{\mathrm{d}}

\newcommand{\bbR}{{\mathbb R}}

\newcommand{\bbP}{{\mathbb P}}

\newcommand{\bbZ}{{\mathbb Z}}

\newcommand{\E}{{\mathrm e}}
\newcommand{\iC}{{\mathrm i}}

\newcommand{\Lie}{\operatorname{\mathsf{L}}}

\newcommand{\la}{\langle}
\newcommand{\ra}{\rangle}

\newcommand{\w}{{\mathchoice{\,{\scriptstyle\wedge}\,}%
{{\scriptstyle\wedge}}{{\scriptscriptstyle\wedge}}%
{{\scriptscriptstyle\wedge}}}}

\newcommand{\be}{\begin{equation}}
\newcommand{\ee}{\end{equation}}
\newcommand{\bpm}{\begin{pmatrix}}
\newcommand{\epm}{\end{pmatrix}}

\numberwithin{equation}{section}

\newtheorem{theorem}{Theorem}

\newtheorem{proposition}{Proposition}
\newtheorem{corollary}{Corollary}

\theoremstyle{remark}
\newtheorem{definition}{Definition}

\newtheorem{remark}{Remark}

\begin{document}

\author[R. Bryant]{Robert L. Bryant}
\address{Duke University Mathematics Department\\
         P.O. Box 90320\\
         Durham, NC 27708-0320}
\email{\href{mailto:bryant@math.duke.edu}{bryant@math.duke.edu}}
\urladdr{\href{http://www.math.duke.edu/~bryant}%
         {http://www.math.duke.edu/\lower3pt\hbox{\symbol{'176}}bryant}}

\title[Geodesically reversible Finsler spheres]
      {Geodesically reversible\\
       Finsler $2$-spheres\\
       of constant curvature}

\date{August 02, 2004}

\begin{abstract}
A Finsler space~$(M,\Sigma)$ is said to be \emph{geodesically
reversible} if each oriented geodesic can be reparametrized
as a geodesic with the reverse orientation.  
A reversible Finsler space is geodesically reversible,
but the converse need not be true.

In this note, building on recent work of LeBrun and Mason~\cite{MR04d:53043},
it is shown that a geodesically reversible Finsler metric
of constant flag curvature on the $2$-sphere is necessarily
projectively flat.

As a corollary, using a previous result of the author~\cite{MR98i:53101}, 
it is shown that a reversible Finsler metric of constant flag 
curvature on the $2$-sphere is necessarily a Riemannian metric
of constant Gauss curvature, thus settling a long-standing problem
in Finsler geometry.
\end{abstract}



\thanks{
Thanks to Duke University for its support via a research grant
and to the NSF for its support via DMS-0103884.
\hfill\break
\hspace*{\parindent} 
This is Version~$2.0$ of ReversibleFinslerSphere.tex.
}

\maketitle

\section{Introduction}\label{sec: intro}

The purpose of this note is to settle a long-standing
problem in Finsler geometry:  Whether there exists a
\emph{reversible} Finsler metric on the $2$-sphere with
constant flag curvature that is not Riemannian.  
By making use of some old results and a fundamental new
result of LeBrun and Mason, I show that such Finsler
structures do not exist.  

First, I prove something related:  
Any geodesically reversible Finsler metric on the $2$-sphere
with constant flag curvature must be projectively flat.
Since the projectively flat Finsler metrics with constant
flag curvature on~$S^2$ were classified some years 
ago~\cite{MR98i:53101}, the above result then reduces
to examining the Finsler structures provided by this
classification.

In a famous 1988 paper~\cite{MR91f:53069}, Akbar-Zadeh 
showed that a (not necessarily reversible) Finsler structure 
on a compact surface with constant negative flag curvature
was necessarily Riemannian or with zero flag curvature 
was necessarily a translation-invariant Finsler structure 
on the standard $2$-torus~$\bbR^2/\bbZ^2$.
This naturally raised the question about what happens
in the case of constant positive flag curvature.

This problem was made more interesting by the discovery
of non-reversible Finsler metrics on the $2$-sphere
with constant positive flag curvature in~\cite{MR97e:53128}.
(However, it should be pointed out that Katok had already
constructed non-reversible Finsler metrics on the 
$2$-sphere~\cite{MR86g:58036} that later turned out to 
have constant flag curvature, although, apparently, 
this was not known at the time of~\cite{MR97e:53128}.)

In the interests of brevity, no attempt has been made to
give an exposition of the basics of Finsler geometry.
There are many sources for this background material however, 
among them \cite{MR01g:53130}, \cite{Ca2}, \cite{MR4259c,MR11212a}, 
and~\cite{MR21:4462}.  

For background more specifically suited for studying
the case of constant flag curvature, including its
proper formulation in higher dimensions, see 
\cite{BS}, \cite{Fu1, Fu2, Fu3}, 
and~\cite{MR97m:53120,MR03i:53111,MR03k:53091}.

The corresponding question about (geodesically) reversible Finsler 
metrics of constant positive flag curvature on the~$n$-sphere 
for~$n>2$ remains open at this writing, since an essential
component of the proof for~$n=2$ that is due to LeBrun and Mason
has not yet been generalized to higher dimensions.

\setcounter{tocdepth}{2}
\tableofcontents

\section{Structure equations}\label{sec: StrEqs}
In this section, Cartan's structure equations for a Finsler
surface will be recalled.

\subsection{Cartan's coframing}
Let~$M$ be a surface and let~$\Sigma\subset TM$ be
a smooth Finsler structure.  I.e., $\Sigma$ 
is a smooth hypersurface in~$M$ such that the basepoint
projection~$\pi:\Sigma\to M$ is a surjective submersion
and such that each fiber
\be
\pi^{-1}(x) = \Sigma_x = \Sigma\cap T_xM
\ee
is a smooth, strictly convex curve in~$T_xM$ whose convex
hull contains the origin~$0_x$ in its interior. 

\begin{remark}[Reversibility]
Note that there is no assumption that~$\Sigma=-\Sigma$.  In other
words, a Finsler structure need not be `reversible' (some 
sources call this property `symmetry'), and assumption
is not needed for the development of the local theory.
\end{remark}

One should think of~$\Sigma$ as the unit vectors 
of a `Finsler metric', i.e., a function~$F:TM\to\bbR$
that restricts to each tangent space~$T_xM$ to be a 
not-necessarily-symmetric but strictly convex Banach norm
on~$T_xM$.  

\subsubsection{$\Sigma$-length of oriented curves}
A curve~$\gamma:(a,b)\to M$ will be said to
be a $\Sigma$-curve (or `unit speed curve') if~$\gamma'(t)$
lies in~$\Sigma$ for all~$t\in(a,b)$.  Any smooth, immersed
curve~$\gamma:(a,b)\to M$ has an orientation-preserving
reparametrization~$h:(u,v)\to(a,b)$ such that~$\gamma\circ h$
is a $\Sigma$-curve.  This reparametrization is unique
up to translation in the domain of~$h$.  Thus, one can 
unambiguously define the (oriented) $\Sigma$-length of a 
subcurve~$\gamma:(\alpha,\beta)\to M$ 
to be~$h^{-1}(\beta)-h^{-1}(\alpha)$, when~$a<\alpha<\beta<b$.

\subsubsection{Cartan's coframing}
The fundamental result about the geometry of Finsler
surfaces is due to Cartan~\cite{Ca1}:

\begin{theorem}[Canonical coframing]\label{thm: cancofr}
Let $\Sigma\subset TM$ be a Finsler structure on the oriented surface~$M$
with basepoint projection~$\pi:\Sigma\to M$.  Then there exists a unique 
coframing~$\bigl(\omega_1,\omega_2,\omega_3\bigr)$ on~$\Sigma$ with the properties:
\begin{enumerate}
\item $\omega_1\w\omega_2$ is a positive multiple of any 
      $\pi$-pullback of a positive $2$-form on~$M$,
\item The tangential lift~$\gamma'$ of any $\Sigma$-curve 
      satisfies~$(\gamma')^*\omega_2=0$ and $(\gamma')^*\omega_1 = dt$,	  
\item $\d\omega_1\w\omega_2=0$,
\item $\omega_1\w \d\omega_1 = \omega_2\w \d\omega_2$, and
\item $\d\omega_1 = \omega_3\w\omega_2$ and $\omega_3\w \d\omega_2 = 0$.
\end{enumerate}
Moreover, there exist unique functions~$I$, $J$, and $K$ on~$\Sigma$
so that
\be\label{eq: CaStrEqs}
\begin{aligned}
\d\omega_1 &= -\omega_2\w\omega_3,\\
\d\omega_2 &= -\omega_3\w\bigl(\omega_1 - I\,\omega_2\bigr),\\
\d\omega_3   &= -\bigl(K\,\omega_1 - J\,\omega_3\bigr)\w\omega_2.\\
\end{aligned}
\ee
\end{theorem}

\begin{remark}[The invariants~$I$, $J$, and~$K$.]
The $1$-form~$\omega_1$ is called {\it Hilbert's invariant integral\/}.
A $\Sigma$-curve~$\gamma$ is a geodesic of the Finsler structure 
if and only if its tangential lift satisfies~$(\gamma')^*\omega_3=0$.  
(Of course, by definition,~$(\gamma')^*\omega_2=0$.)

The function~$I$ vanishes if and only if~$\Sigma$ is the unit circle bundle
of a Riemannian metric on~$M$, in which case the function~$K$ becomes
the $\pi$-pullback of the Gauss curvature of the underlying metric.  

The function~$J$ vanishes if and only if the Finsler structure 
is what is called {\it Landsberg\/}~\cite{MR01g:53130}.

The function~$K$ is known as the Finsler-Gauss curvature and plays
the same role in the Jacobi equation for Finsler geodesics as the
Gauss curvature does in the Jacobi equation for Riemannian geodesics.
\end{remark}

Let~$X_1$, $X_2$, and~$X_3$ be the vector fields on~$\Sigma$ that
are dual to the coframing~$(\omega_1,\omega_2,\omega_3)$.  
Then the flow of~$X_1$ is the geodesic flow on~$\Sigma$.

\begin{remark}[The effect of orientations]\label{rem: orientation}
If one reverses the orientation of~$M$, then the canonical
coframing~$\omega$ on~$\Sigma$ 
is replaced by~$(\omega_1,-\omega_2,-\omega_3)$.  

In fact, Cartan's actual statement of Theorem~\ref{thm: cancofr}
does not assume that~$M$ is oriented and concludes that 
there is a canonical coframing on~$\Sigma$ up to the sign 
ambiguity given here. The present version of the statement
is a trivial rearrangement of Cartan's that is more easily
applied in the situations encountered in this note.
\end{remark}

\subsubsection{Reconstruction of~$M$ and its Finsler structure}
The information contained in the $3$-manifold~$\Sigma$ and
its coframing~$\omega = (\omega_1,\omega_2,\omega_3)$ is sufficient to
recover~$M$, its orientation, and the embedding of~$\Sigma$ into~$M$, 
a fact that is implicit in Cartan's analysis:

\begin{proposition}[Isometries and automorphisms]
\label{prop: Iso&Auto}
For any orientation-preserving Finsler isometry~$\phi:M\to M$,
its derivative~$\phi':TM\to TM$ induces 
a diffeomorphism~$\phi':\Sigma\to\Sigma$ 
that preserves the coframing~$\omega  = (\omega_1,\omega_2,\omega_3)$.

Conversely, any diffeomorphism~$\psi:\Sigma\to\Sigma$ that
preserves~$\omega$ is of the form~$\psi=\phi'$ for a unique
orientation-preserving Finsler isometry~$\phi:M\to M$.
\end{proposition}

\begin{proof}
The first statement follows directly from~Theorem~\ref{thm: cancofr}.
I will sketch how the converse goes.

The integral curves of the system~$\omega_1=\omega_2=0$
on~$\Sigma$ are closed and the codimension~$2$ foliation they define 
has trivial holonomy, so $M$ can be identified with the leaf
space of this system and carries a unique smooth structure for
which the leaf projection~$\pi:\Sigma\to M$ is a smooth submersion.

Because of the connectedness of the $\pi$-fibers, there will be a unique
orientation on~$M$ such that a positive $2$-form pulls back under~$\pi$
to be a positive multiple of~$\omega_1\w\omega_2$.  Thus, $M$,
its smooth structure, and its orientation can be recovered from
the coframing.

The inclusion~$\iota:\Sigma\to TM$ is then seen to be simply 
given by~$\iota(u) = \pi'\bigl(X_1(u)\bigr)\in T_{\pi(u)}M$.
Thus, even the Finsler structure on~$M$ can be recovered from~$\Sigma$
and the coframing.

The desired result now follows by noting that any~$\psi:\Sigma\to\Sigma$
that preserves~$\omega$ will necessarily preserve the integral curves
of the system~$\omega_1=\omega_2=0$ and hence induce a map~$\phi:M\to M$
that is $\pi$-intertwined with~$\psi$.  The verification that~$\phi$
is an orientation-preserving Finsler isometry is easy and can be left
to the reader.
\end{proof}

\begin{corollary}[Orientation-reversing isometries]
\label{cor: orrevisoms}
Any diffeomorphism $\psi:\Sigma\to\Sigma$ that
satisfies~$\psi^*(\omega) = (\omega_1,-\omega_2,-\omega_3)$
is of the form~$\psi=\phi'$ for a unique
orientation-\emph{reversing} Finsler isometry~$\phi:M\to M$.
\end{corollary}

\subsection{Bianchi identities}
Taking the exterior derivatives 
of the structure equations~\eqref{eq: CaStrEqs}
yields the formulae
\be\label{eq: BianchiIds}
\bpm \d I \\ \d J\\ \d K\epm
= \bpm J & I_2 & I_3 \\
      -K_3-KI & J_2 & J_3 \\
      K_1 & K_2 & K_3  \epm
\bpm \omega_1 \\ \omega_2 \\ \omega_3\epm
\ee
for some functions~$I_2$, $I_3$, $J_2$, $J_3$, $K_1$, $K_2$, and~$K_3$
on~$\Sigma$.

\subsection{Simplifications when $K\equiv1$}
The Finsler structures of interest in this article are the 
ones that satisfy~$K\equiv1$.  In this case, the structure
equations simplify to
\be\label{eq: CaStrEqsK=1}
\begin{aligned}
\d\omega_1 &= -\omega_2\w\omega_3,\\
\d\omega_2 &= -\omega_3\w\bigl(\omega_1 - I\,\omega_2\bigr),\\
\d\omega_3 &= -\bigl(\omega_1 - J\,\omega_3\bigr)\w\omega_2,\\
\end{aligned}
\ee
and the Bianchi identities become
\be\label{eq: BianchiIdsK=1}
\bpm \d I \\ \d J\epm
= \bpm J & I_2 & I_3 \\
      -I & J_2 & J_3\epm
\bpm \omega_1 \\ \omega_2 \\ \omega_3\epm.
\ee

\begin{remark}[A geodesic conservation law]
The equations~\eqref{eq: BianchiIdsK=1} imply that
the function $I^2+J^2$ is constant on the integral
curves of~$\omega_2=\omega_3=0$, i.e., the lifts of geodesics. 
This function need not be constant on~$\Sigma$, 
in which case, it provides a nontrivial conservation law 
for the geodesic flow on~$\Sigma$.  (Of course, this function
vanishes identically in the Riemannian case.)
\end{remark}

\subsection{Some global consequences of~$K\equiv1$}
Suppose now that~$M$ is connected and geodesically complete, 
i.e., that, the vector field~$X_1$ is complete on~$\Sigma$
(in both forward and backward time).  Of course, if~$M$
were assumed to be compact, then~$\Sigma$ would be 
also, and the completeness of~$X_1$ would follow from this.

The assumption that~$M$ be connected implies
that~$\Sigma$ is connected.

Let~$\Psi:\Sigma\times\bbR\to\Sigma$ be the flow of~$X_1$
and, for brevity, let~$\Psi_t:\Sigma\to\Sigma$ denote 
the time~$t$ flow of~$X_1$. 
Since the structure equations imply
\be
\Lie_{X_1}\omega_1=0,\qquad
\Lie_{X_1}\omega_2=\omega_3,\qquad
\Lie_{X_1}\omega_3=-\omega_2,\qquad
\ee 
it follows (letting~$t:\Sigma\times\bbR\to\bbR$ denote
the coordinate that is the projection on the second factor) that
\be\label{eq: Psi_tpullbacks}
\begin{aligned}
\Psi^*\omega_1 & = \omega_1 + \d t, \\
\Psi^*\omega_2 & = \phantom{-}\cos t\,\omega_2 + \sin t\,\omega_3, \\
\Psi^*\omega_3 & = -\sin t\,\omega_2 + \cos t\,\omega_3. \\
\end{aligned}
\ee

\begin{proposition}[The quasi-antipodal map]
\label{prop: quasiantipodalmap}
There exists a unique orientation-reversing Finsler
isometry~$\alpha:M\to M$ such that~$\alpha' = \Psi_\pi$. 
For any point~$p\in M$, every unit speed geodesic leaving~$p$
passes through~$\alpha(p)$ at distance~$\pi$.
\end{proposition}

\begin{proof}
By~\eqref{eq: Psi_tpullbacks}, it follows that
$\Psi_\pi:\Sigma\to\Sigma$ satisfies
\be
\Psi^*_\pi\omega = (\omega_1,-\omega_2,-\omega_3).
\ee
Hence, by Corollary~\ref{cor: orrevisoms}, 
there is a unique orientation-reversing Finsler isometry~$\alpha:M\to M$
such that $\Psi_\pi = \alpha':\Sigma\to\Sigma$.

Since~$X_1$ is the geodesic flow vector field, 
any unit speed geodesic leaving~$p$ at time~$0$ is
of the form~$\gamma(t) = \pi\bigl(\Psi_t(u)\bigr)$
for some~$u\in\Sigma_p\subset T_pM$.  Thus,~$\gamma(\pi)
= \pi\bigl(\Psi_\pi(u)\bigr)=\pi\bigl(\alpha'(u)\bigr)=\alpha(p)$,
as claimed.
\end{proof}

Now, for any fixed~$p\in M$, the fiber~$\Sigma_p\subset T_pM$, 
is diffeomorphic to a circle and is naturally oriented
by taking the pullback of~$\omega_3$ to~$\Sigma_p$ to be
a positive $1$-form.   Define~$r(p)>0$ by
\be
r(p) = \frac1{2\pi}\int_{\Sigma_p} \omega_3\,.
\ee
Then~$\Sigma_p$ can be parametrized 
by a mapping~$\iota_p:S^1\to\Sigma_p$ 
that satisfies~$\iota_p^*(\omega_3) =r(p)\,\d\theta$ and
that is uniquely determined once one fixes~$\iota_p(0)=u\in \Sigma_p$.
Such a parametrization~$\iota_p$ will be referred to as
an \emph{angle measure} on~$\Sigma_p$.

\begin{proposition}[Geodesic polar coordinates]
\label{prop: geodpolarcoords}
For any~$p\in M$, fix an angle measure~$\iota_p:S^1\to\Sigma_p$.
Then the mapping~$E_p:S^2\to M$ defined by
\be\label{eq: Epdefined}
E_p(\sin t\cos\theta,\sin t\sin\theta,\cos t)
= \pi\bigl(\Psi_t(\iota_p(\theta))\bigr)
\ee
is an orientation-preserving homeomorphism 
that is smooth away from~$(0,0,\pm1)\in S^2$.  
In particular, $M$ is homeomorphic to the $2$-sphere 
and its diameter as a Finsler space is equal to~$\pi$.
\end{proposition}

\begin{proof}
Consider the mapping~$R_p:S^1\times\bbR\to\Sigma$ defined by
\be
R_p(\theta,t) = \Psi\bigl(\iota_p(\theta),t\bigr).
\ee
The formulae~\eqref{eq: Psi_tpullbacks},
the fact that~$\Psi$ is the flow of~$X_1$, and the
defining property of~$\iota_p$ then combine to show that
\be
R_p^*(\omega_1\w\omega_2) 
= \d t \w\bigl(\sin t\,\,r(p)\,\d\theta\bigr)
= r(p)\,\sin t\,\d t\w\d\theta.
\ee
Thus, the composition 
$\pi{\circ}R_p: S^1\times\bbR\to M$ is a smooth map that
is a local diffeomorphism
away from the circles~$(\theta,t) = (\theta,k\pi)$ for each integer~$k$.
Of course,~$\pi\bigl(R_p(\theta,0)\bigr)=p$ and 
$\pi\bigl(R_p(\theta,\pi)\bigr)=\alpha(p)$ for all~$\theta\in S^1$.

It now follows that the formula~\eqref{eq: Epdefined} 
well-defines a mapping~$E_p:S^2\to M$ that is smooth 
and an orientation-preserving local diffeomorphism 
away from~$(0,0,\pm1)$.  Near the two points~$(0,0,\pm1)$, 
the mapping~$E_p$ is still a (not necessarily differentiable) 
orientation-preserving local homeomorphism.

It follows that~$E_p:S^2\to M$ is a topological covering map.
Since~$M$ is orientable by assumption, it follows that~$E_p$
must be a homeomorphism and, in particular, must be one-to-one
and onto.  The statement about diameters follows.
\end{proof} 

\begin{remark}
Versions of Propositions~\ref{prop: quasiantipodalmap} 
and~\ref{prop: geodpolarcoords} 
were proved by Shen~\cite{MR97m:53120}
in the case that~$\Sigma$ is reversible 
(see Definition~\ref{def: reversibility}).
\end{remark}

\begin{proposition}\label{prop: alphaproperties}
Either~$\alpha^2=\mathrm{id}$ on~$M$ {\upshape(}in which case, all of the
$\Sigma$-geodesics are closed of length~$2\pi${\upshape)} or else~$\alpha^2$
has exactly two fixed points, say~$n$ and~$\alpha(n)$.

In the latter case, there exists a positive definite inner product
on~$T_nM$ that is invariant under~$(\alpha^2)'(n):T_nM\to T_nM$ and
there is an angle~$\theta_n\in(0,2\pi)$ such that~$(\alpha^2)'(n)$
is a counterclockwise rotation by~$\theta_n$ in this inner product.
\end{proposition}

\begin{proof}
Assume that~$\alpha^2:M\to M$ is not the identity, or else there
is nothing to prove.  Since~$\alpha^2$ is an orientation preserving
diffeomorphism of the $2$-sphere, it must have at least one fixed
point.  Let~$n$ be such a fixed point.  By the very definition 
of~$\alpha$, it then follows that~$\alpha(n)$ is also a fixed point
of~$\alpha^2$.  It must be shown that $\alpha^2$ 
has no other fixed points.

First, consider the linear map~$L = (\alpha^2)'(n):T_nM\to T_nM$.
Since~$\alpha^2$ is a Finsler isometry, 
the linear map~$L$ must preserve~$\Sigma_n\subset T_nM$.  
Let~$K_n\subset T_nM$ be the convex set bounded by~$\Sigma_n$.

Define a positive definite quadratic form on~$T^*_nM$ 
by letting~$\la\lambda_1,\lambda_2\ra$ be
defined for~$\lambda_1,\lambda_2\in T^*_nM$ 
to be the average of the quadratic function~$\lambda_1\lambda_2$ 
over~$K_n$ (using any translation 
invariant measure on~$K_n$ induced by its inclusion into the vector 
space~$T_nM$).  Since~$L$ is a linear map carrying~$K_n$ into itself,
it must preserve this quadratic form and hence must also
preserve the dual (positive-definite) quadratic form on~$T_nM$.  
Since~$L$ also preserves an orientation on~$T_nM$, 
it follows that, with respect to this invariant inner product,
$L$ must be a counterclockwise rotation by
some angle~$\theta_n\in[0,2\pi)$.  

If~$\theta_n$ were~$0$, i.e.,~$L$ were the identity on~$T_nM$, 
then all of the geodesics through~$n$ would close at length~$2\pi$.
In particular, the mapping~$\Psi_{2\pi}:\Sigma\to\Sigma$ would
have a fixed point and would preserve the coframing~$\omega$,
implying that~$\Psi_{2\pi}$ is the identity on~$\Sigma$ 
and hence that~$\alpha^2$ would be the identity.  
Thus,~$0<\theta_n<2\pi$.

Since~$n$ was an arbitrarily chosen fixed point of~$\alpha^2$,
it follows that every fixed point of~$\alpha^2$ is an isolated
elliptic fixed point, i.e., a fixed point of index~$1$.  
Since~$M$ is diffeomorphic to~$S^2$, 
the Hopf Index Theorem implies that
the map~$\alpha^2$ has exactly two fixed points.  
Thus~$\alpha^2$ has no fixed points 
other than~$n$ and~$\alpha(n)$.
\end{proof}

\begin{remark}[The Katok examples]
The Katok examples analyzed by Ziller~\cite{MR86g:58036} turn out%
\footnote{Colleen Robles, private communication}
to have $K\equiv1$ and are examples in which~$\alpha^2$ 
is not the identity.  Thus, the second possibility
in Proposition~\ref{prop: alphaproperties} does occur.

In any case, when~$\alpha^2$ is not the identity,
~$\theta_n+\theta_{\alpha(n)} = 2\pi$.

If the angle~$\theta_n$ defined in Proposition~\ref{prop: alphaproperties}
is not a rational multiple of~$\pi$, then the iterates of~$\alpha^2$
are dense in a circle of Finsler isometries of~$(M,\Sigma)$ that
fix~$n$ and~$\alpha(n)$.  In such a case, $(M,\Sigma)$ is rotationally
symmetric about~$n$.  Moreover, it is symmetric (in an orientation
reversing sense) with respect to~$\alpha$.

If~$\theta_n = 2\pi(p/q)$ where~$0<p\le q$ and $p$ and~$q$ have
no common factors, then $\alpha^{2q}$ is the identity, so that every
geodesic closes at length~$2\pi q$ (though some may close sooner).
\end{remark}

\section{A double fibration}\label{sec: DblFibr}

Throughout this section~$\Sigma$ will be assumed to be
a Finsler structure on~$M$ (assumed diffeomorphic to the $2$-sphere)
satisfying~$K\equiv1$.

I begin by noting that, if all the geodesics on~$M$ close at
distance~$2\pi$, then the set of oriented $\Sigma$-geodesics 
has the structure of a manifold in a natural way.

\begin{proposition}[The space of oriented geodesics]
\label{prop: Lambdadefined}
If~$\alpha^2$ is the identity, 
then the action
\be\label{eq: S1action}
u\cdot \E^{\iC t } = \Psi(u,t)
\ee
defines a smooth, free $S^1$-action on~$\Sigma$ whose
orbits are the integral curves of~$\omega_2=\omega_3=0$ and
there exists a smooth surface~$\Lambda$ diffeomorphic to~$S^2$
and a smooth submersion~$\lambda:\Sigma\to\Lambda$
so that the action~\eqref{eq: S1action} makes 
$\lambda:\Sigma\to\Lambda$ into a principal right~$S^1$-bundle
over~$\Lambda$.
\end{proposition}

\begin{proof}
If~$\alpha^2$ is the identity, then the flow of~$X_1$ is
periodic of period~$2\pi$, so~\eqref{eq: S1action} defines
a smooth $S^1$-action on~$\Sigma$.  Since~$X_1$ never
vanishes, this action has no fixed points.  Thus, if this action 
were not free, then there would be a~$u\in\Sigma$ and an 
integer~$k\ge2$ such that~$\Psi(u,2\pi/k)=u$.  However, 
since~$0<2\pi/k\le \pi$, the equality~$\Psi(u,2\pi/k)=u$
would violate Proposition~\ref{prop: geodpolarcoords}, since
then~$E_{\pi(u)}:S^2\to M$ could not be one-to-one.  

Thus, the $S^1$-action~\eqref{eq: S1action} is free 
and the rest of the proposition follows by standard arguments.
\end{proof}

\begin{remark}[Double fibration and path geometries]
The two mappings~$\pi:\Sigma\to M$ and~$\lambda:\Sigma\to\Lambda$
define a double fibration and it is easy to see that this 
double fibration satisfies the usual nondegeneracy axioms for
double fibrations. For example,~$\lambda\times\pi:\Sigma\to \Lambda\times M$ 
is clearly a smooth embedding.  The other properties are similarly
easy to verify using the structure equations.  Thus,~$\Sigma$ defines
a (generalized) path geometry on each of~$\Lambda$ and~$M$.

For more background on path geometries and their invariants,
see, for example, Section~2 of~\cite{MR98i:53101}.
\end{remark}

\subsection{Induced structures on~$\Lambda$}
I will now recall some results from~\cite{MR98i:53101}.
Throughout this subsection, I will be assuming that~$\alpha^2$
is the identity, so that~$\Lambda$ exists as a smooth manifold.

The relations~\eqref{eq: Psi_tpullbacks} 
show that the quadratic form~${\omega_2}^2+{\omega_3}^2$ 
is invariant under the flow of~$X_1$.  
Consequently, there is a unique Riemannian metric
on~$\Lambda$, say~$g$, such that
\be
\lambda^*(g) = {\omega_2}^2+{\omega_3}^2
\ee
Moreover, the $2$-form~$\omega_3\w\omega_2$ is invariant
under the flow of~$X_1$, so it is the pullback under~$\lambda$
of an area $2$-form for~$g$, which will be denoted~$\d A_g$. 

Now, there is an embedding~$\xi:\Sigma\to T\Lambda$ defined by
\be
\xi(u) = \lambda'\bigl(X_3(u)\bigr)
\ee
and one sees that~$\xi$ embeds~$\Lambda$ as the unit
sphere bundle of~$\Lambda$ endowed with the metric~$g$.

The structure equations~\eqref{eq: CaStrEqsK=1} show that, 
under this identification of~$\Sigma$ with the unit
sphere bundle of~$\Lambda$, the Levi-Civita connection form 
on~$\Sigma$ is
\be
\rho = -\omega_1 + I\,\omega_2 + J\,\omega_3\,.
\ee
Note that~$-\omega_1$ and $I\,\omega_2 + J\,\omega_3$
are invariant under the flow of~$X_1$.  

For the next two results, which follow from the 
structure equations derived so far by simply
unraveling the definitions, the reader may want to consult
LeBrun and Mason~\cite{MR04d:53043} for the definition
and properties of the projective structure 
associated to an affine connection on a surface.
[They restrict themselves to the consideration of
torsion-free connections, but, as they point out, 
this does not affect the results.]

\begin{proposition}\label{prop: nablaonLambda}
There exists a $g$-compatible affine connection~$\nabla$ 
on~$\Lambda$ such that the $\nabla$-geodesics 
are the $\lambda$-projections of the integral curves 
of~$\omega_1=\omega_2=0$.
\qed
\end{proposition}

\begin{corollary}
The geodesics of the projective structure~$[\nabla]$
on~$\Lambda$ are closed.
\end{corollary}

\begin{proof}
By Proposition~\ref{prop: nablaonLambda}, the geodesics 
of~$[\nabla]$ are the $\lambda$-projections of the integral
curves of the system~$\omega_1=\omega_2=0$, but these
integral curves are the fibers of the map~$\pi:\Sigma\to M$ 
and hence are closed.
\end{proof}

\subsection{Geodesic reversibility implies geodesic periodicity}
It is now time to come to the main point of this note.

\begin{definition}[Reversibility]\label{def: reversibility}
The Finsler structure~$\Sigma\subset TM$
is said to be \emph{reversible} if~$\Sigma = - \Sigma$.
\end{definition}

\begin{definition}[Geodesic reversibility]
A Finsler structure~$\Sigma\subset TM$ will be said to
be \emph{geodesically reversible} if any $\Sigma$-geodesic
$\gamma:(a,b)\to TM$ can be reparametrized in an 
orientation-reversing way so as to remain a $\Sigma$-geodesic. 
\end{definition}

\begin{remark}
Any reversible Finsler structure is geodesically
reversible.  On the other hand,
the non-Riemannian Finsler examples constructed in Section~4
of~\cite{MR98i:53101} are geodesically reversible 
but not reversible, so the reverse implication does not hold.
\end{remark}

\begin{proposition}\label{prop: geodrevimpliesclosure}
If~$(M,\Sigma)$ is geodesically reversible, then~$\alpha^2$ 
is the identity on~$M$.
\end{proposition} 

\begin{proof}
For any point~$p\in M$, consider the geodesics leaving~$p$.
By Proposition~\ref{prop: geodpolarcoords}, 
they all converge at distance~$\pi$
on~$\alpha(p)$ but do not intersect between distance~$0$ 
and distance~$\pi$.  By assumption, reversing these geodesic
segements, i.e., tracing them backwards from~$\alpha(p)$, 
yields $\Sigma$-geodesics (which are no longer necessarily 
unit speed).  Moreover, all of these 
geodesics remain disjoint until they pass through~$p$, 
at which point, they all converge.

However, again by Proposition~\ref{prop: geodpolarcoords}, 
the unit speed geodesics leaving~$\alpha(p)$ remain disjoint 
for distances between~$0$ and~$\pi$ and they all converge 
on~$\alpha\bigl(\alpha(p)\bigr)$ at distance~$\pi$.

It follows that~$\alpha\bigl(\alpha(p)\bigr)$ must be~$p$.
In other words,~$\alpha^2$ is the identity.
\end{proof}

\begin{remark}
The converse of Proposition~\ref{prop: geodrevimpliesclosure} 
does not hold.  The~$K\equiv1$ examples provided by 
Theorem~$3$ of~\cite{bry3} that are based on Guillemin's Zoll 
metrics have all their geodesics closed of length~$2\pi$ 
(and hence $\alpha^2$ is the identity), 
but none of the non-Riemannian ones are geodesically reversible. 
\end{remark}

\subsection{Geodesic reversibility implies projective flatness}
The next step is to consider the space of \emph{unoriented}
$\Sigma$-geodesics on~$M$.  This only makes sense if one assumes
that~$\Sigma$ is geodesically reversible, so assume this for
the rest of this subsection.

For each oriented $\Sigma$-geodesic~$\gamma:S^1\to M$, let
$\beta(\gamma)$ denote the reversed curve, reparametrized so
as to be a~$\Sigma$-geodesic.  Obviously~$\beta:\Lambda\to\Lambda$
is a fixed-point free involution of~$\Lambda$, so that the
quotient manifold~$\Lambda/\beta$ is diffeomorphic to~$\bbR\bbP^2$.

\begin{proposition}\label{prop: projstronLmodb}
The path geometry on~$\Lambda$ defined by the geodesics of~$[\nabla]$
is invariant under~$\beta$ and hence descends to a well-defined
path geometry on~$\Lambda/\beta$.  Moreover, this path geometry
is the path geometry of a projective connection on~$\Lambda/\beta$
with all of its geodesics closed.
\end{proposition}

\begin{proof} 
Since, by definition, a point~$p$ in~$M$ lies on a geodesic~$\gamma$
if and only if it lies on~$\beta(\gamma)$, it follows that~$\beta$
carries each~$[\nabla]$-geodesic into itself.  In particular, even
though~$\beta$ may not (indeed, most likely does not) preserve~$\nabla$, 
it must preserve~$[\nabla]$ since the projective equivalence class
of~$\nabla$ is determined by its geodesics.  Thus, the claims of
the Proposition are verified.
\end{proof}

It is at this point that the crucial contribution of LeBrun and
Mason~\cite{MR04d:53043} enters:

\begin{theorem}[LeBrun-Mason]
\label{thm: LeBrunMason}
Any projective structure on~$\bbR\bbP^2$ that has all of its
geodesics closed is projectively equivalent to the standard
{\upshape(}i.e., flat{\upshape)} projective structure.
\end{theorem}

\begin{corollary}\label{cor: Lambdaprojflat}
If~$\Sigma$ is a geodesically reversible Finsler structure
on~$M\simeq S^2$ that satisfies~$K\equiv1$, then the
induced projective structure~$[\nabla]$ on~$\Lambda$ is
projectively flat. \qed
\end{corollary}

\begin{remark}[LeBrun and Mason's classification]
The article~\cite{MR04d:53043} contains, 
in addition to Theorem~\ref{thm: LeBrunMason}, 
much information about \emph{Zoll projective structures}
on the $2$-sphere, i.e., projective structures on the $2$-sphere
all of whose geodesics are closed.  It turns out that, 
in a certain sense, there are many more of them than 
there are Zoll metrics on the $2$-sphere.  

Their results could quite likely be very useful 
in understanding the case of non-reversible Finsler
metrics satisfying~$K\equiv1$ on the $2$-sphere
that satisfy~$\alpha^2=\mathrm{id}$, which is still
not very well understood.  It is even possible that an
orbifold version of their results could be useful in 
the case in which~$\alpha^2$ is not the identity but
has finite order.  This may be the subject of a later
article.
\end{remark}

\section{Classification}\label{sec: Classification}

In this final section, the main theorem will be proved.

\subsection{Consequences of projective flatness}
Recall from Section~$2$ of~\cite{MR98i:53101} that
if a projective structure on a surface is projectively
flat then its dual path geometry is projective and,
moreover, projectively flat.

\begin{proposition}\label{prop: MSigmaprojflat}
If~$\Sigma$ is a geodesically reversible Finsler structure
on~$M\simeq S^2$ with~$K\equiv1$, then the $\Sigma$-geodesics
in~$M$ are the geodesics of a flat projective structure.
\end{proposition}

\begin{proof}
The dual path geometry of~$\Lambda$ with its projective
structure~$[\nabla]$ is $M$ with the space of paths being
the $\Sigma$-geodesics.  Now apply Corollary~\ref{cor: Lambdaprojflat}.
\end{proof}

\begin{corollary}\label{cor: geodreversclassify}
Let~$M$ be diffeomorphic to~$S^2$.
Up to diffeomorphism, any geodesically reversible Finsler
structure~$\Sigma\subset TM$ with~$K\equiv 1$ is equivalent
to a member of the $2$-parameter family described in Theorem~$10$
of~\cite{MR98i:53101}.
\end{corollary}

\begin{proof}
In light of Proposition~\ref{prop: MSigmaprojflat}, one can
apply Theorems~$9$ and~$10$ of~\cite{MR98i:53101}, which
gives the result.
\end{proof}

\begin{remark}
It is interesting to note that each member of the 
$2$-parameter family described in Theorem~$10$
of~\cite{MR98i:53101} is projectively flat and
hence geodesically reversible.
\end{remark}

\subsection{Reversibility}
Now for the main rigidity theorem.

\begin{theorem}\label{thm: reversibleFinslerrigid}
Any reversible Finsler structure
on~$M\simeq S^2$ that satisfies~$K\equiv1$
is Riemannian and hence isometric to the standard unit sphere.
\end{theorem}

\begin{proof}
Such a Finsler structure would be geodesically reversible
and hence, by Corollary~\ref{cor: geodreversclassify}, 
a member of the family described in Theorem~$10$
of~\cite{MR98i:53101}.  However, by inspection, the only
member of this geodesically reversible family that is
actually reversible is the Riemannian one.
\end{proof}

\begin{remark}[The argument of Foulon-Reissman]
In Section~$4$ of~\cite{Foulon2002}, P. Foulon sketches an argument,
due to himself and A. Reissman, that a reversible
Finsler metric on the $2$-sphere satisfying~$K\equiv1$ 
that satisfies a certain integral-geometric condition 
(called by them `Radon-Gelfand') is necessarily Riemannian.
Their condition holds, in particular, whenever the projective
structure~$[\nabla]$ on~$\Lambda$ is projectively flat.
Thus, an alternate proof of Theorem~\ref{thm: reversibleFinslerrigid}
could be given by combining LeBrun and Mason's Theorem~\ref{thm: LeBrunMason}
with Foulon and Reissman's argument.

The proof of Theorem~\ref{thm: reversibleFinslerrigid} 
in this article instead relies on the classification 
in~\cite{MR98i:53101}.  
\end{remark}

\bibliographystyle{hamsplain}

\providecommand{\bysame}{\leavevmode\hbox to3em{\hrulefill}\thinspace}

\end{document}